\documentclass[12pt]{article}
\usepackage{amsfonts,amsmath,amsthm,amscd}

\renewcommand{\a}{\alpha}
\renewcommand{\L}{\Lambda}
\newcommand{\f}{\varphi}
\newcommand{\s}{\sigma}

\newcommand{\G}{\Gamma}

\newcommand{\F}{\mathbb{F}}
\newcommand{\K}{\mathbb{K}}
\newcommand{\Z}{\mathbb{Z}}
\renewcommand{\P}{\mathbb{P}}
\newcommand{\op}{\oplus}
\newcommand{\MZ}{\mathcal Z}
\newcommand{\Gal}{\mathrm{Gal}}
\newcommand{\PhiKZ}{\Phi_{KZ}}
\theoremstyle{definition}
\newtheorem*{thm}{Theorem}
\newtheorem{lemma}{Lemma}
\DeclareMathOperator{\indeg}{indeg}
\DeclareMathOperator{\outdeg}{outdeg}

\begin{document}

\title{Automorphisms of necklaces and sandpile groups}

\author{S.\,Duzhin\and D.\,Pasechnik
\thanks{Supported by Singapore MOE Tier 2 Grant MOE2011-T2-1-090,
the first author also by grant NSh-4850.2012.1.}}

\date{}

\maketitle

\begin{abstract}
We introduce a group naturally acting on aperiodic necklaces of length $n$
with two
colours using the 1--1 correspondences between aperiodic necklaces 
and irreducible polynomials over the field $\F_2$ of two elements. 
We notice that this 
group is isomorphic to the quotient group of non-degenerate
circulant matrices of size $n$ over that 
field modulo a natural cyclic subgroup. Our groups turn out to be isomorphic 
to the sandpile groups for a special sequence of directed graphs.
\end{abstract}

\section{Introduction}
\label{intr}

This work originated in the research of the first author related to the
Drinfeld associator \cite{DD}.  It is well-known (see, e.g.  \cite{CDM})
that the logarithm of the classical associator $\PhiKZ$ is an element of the
completed free Lie algebra on two generators with coefficients in the
algebra $\MZ$ of multiple zeta values \cite{Hoff}.  The free Lie algebra has
a basis whose elements are labelled by aperiodic necklaces \cite{Reu}. 
However, the explicit expansion of $\log\PhiKZ$ over this basis shown in
\cite{DD} displays a highly chaotic behaviour.  Therefore, a natural idea
arises to introduce some structure in the set of necklaces.  In this paper,
we introduce a group acting on the set of aperiodic necklaces of fixed length
in the hope that the orbits of this group may shed some light on the
structure of the embarrassing expression given in \cite{DD}.

In the spirit proclaimed by V.\,I.\,Arnold in \cite{Arn}, this paper does not contain any proofs, only constructions, problems, 
motivations, examples
and statement of results proven elsewhere. It may be considered as an 
informal introduction to the paper \cite{CHP}.

\section{General setting}
\label{gen}

In the most general setting, the approach we follow here consists in the
following.

Suppose we have two sets $A$ and $B$ and a family of bijections
$\f_i:A\to B$. This data can be used to define two groups: $G_A$, acting 
on the set $A$, and $G_B$, acting on the set $B$. Indeed, let $G_A$ be the
group generated by all bijections $\f_{ij}=\f_i^{-1}\circ\f_j$ and $G_B$ the
group generated by all $\f_{ij}'=\f_i\circ\f_j^{-1}$. 

The following lemma is immediate.
\begin{lemma}
The groups $G_A$ and $G_B$ are isomorphic, an isomorphism being given by the
assignment $g\mapsto \f_i\circ g\circ\f_i^{-1}$ for any fixed $i$.  In
particular, the actions of these groups on the sets $A$ and $B$ are 
equivalent,
any bijection $\f_i$ maps the orbits of the group $G_A$ onto the 
orbits
of the group $G_B$, and this map on the set of orbits does not depend on the
choice of a particular $i$.
\end{lemma}

\section{History}
\label{hist}

The story began when the first
author obtained the results of Section \ref{reut} (and also 
Section \ref{gol} which is, however, irrelevant to the main topic of 
this paper). This was done experimentally, by computer. The second 
author found the sequence of orders
of the groups $RG_n^2$ in Sloane's encyclopaedia (A027362) and conjectured
that the groups are isomorphic to the sandpile groups of generalized de
Bruijn graphs. This conjecture was checked by computer up to order 16 by 
the present authors and then proved by S.H.Chan 
in his Bachelor's thesis \cite{Ch}.\footnote{Actually, the correct 
proof appeared only
in the updated version of the thesis paper.} Finally, this theorem 
was generalized and extended in the paper \cite{CHP}.

\section{Groups acting on necklaces}
\label{grp}

A $p$-coloured \textit{necklace} of length $n$ is a sequence of $n$ objects of $p$
different kinds (called beads) considered up to cyclic shifts.
Necklaces may be periodic (admitting a non-trivial shift that does not
change it) or aperiodic.
If the colours are linearly ordered, then the lexicographically minimal
of all cyclic shifts of an aperiodic necklace is called a \textit{Lyndon} word. 
We will always suppose that
$p$ is prime and use the elements of a prime finite field $\F_p$ as the
colours of beads. We will refer to such necklaces and Lyndon words as
\textit{arithmetical}.

The main case which is most interesting from the point of view of 
Section \ref{intr} is $p=2$.

Example:
There are exactly 6 aperiodic 2-coloured necklaces of length 5 given by
Lyndon words 00001, 00011, 00111, 01111, 00101, 01011. 

Denote the set of aperiodic necklaces with parameters $n$ and $p$ by $N_n^p$.
There is a formula for the cardinality of this set in terms of the 
M\"{o}bius function:
$$
  |N_n^p|=\frac{1}{n}\sum_{d|n}\mu(d)p^{n/d}.
$$
It is remarkable that this number
is equal to the number of irreducible polynomials
over the field $\F_p$ of degree $n$, or, which is the same, to the 
number of
orbits of maximal length $n$ of the action of the Galois group
$\Gal(\F_{p^n}/\F_p)$. \footnote{This is why we consider only the necklaces with a prime number of colors. Generalizing the theory for $p^k$ instead of $p$ is an interesting open problem.}
Denote the set of irreducible polynomials by 
$I_n^p$. There are two known
explicit constructions of a 1--1 correspondence
$$
  \f: N_n^p \to I_n^p
$$
One of them (first mentioned in Reutenauer \cite{Reu} who ascribed it to
folklore)
depends on the choice of a \textit{normal} polynomial of degree $n$ over
$\F_p$. Another one belongs to Golomb \cite{Gol}
and depends on the choice of a \textit{primitive} polynomial
of degree $n$ over $\F_p$. According to Section \ref{gen}, either set of 
bijections generates a group of transformations on the set of aperiodic 
necklaces. We will call the first one the \textit{Reutenauer} group 
and denote it by $RG_n^p$, and the second, the \textit{Golomb} group 
and denote it by $GG_n^p$.

\section{Galois group and irreducible polynomials}

It is well known that the Galois group $\Gal(\F_{p^n}:\F_p)$ is a cyclic
group of order $n$ generated by the Frobenius automorphism $\s(x)=x^p$.
Each orbit of this group coincides with the set of roots of an irreducible
polynomial over $\F_p$ of degree which divides $n$. The orbits of maximum
length $n$ correspond to irreducible polynomials of degree exactly $n$.
The union of all such orbits is equal to the complement in $\F_{p^n}$
of all its proper subfields. A polynomial is called \textit{normal}, if the
set of its roots constitutes a basis of the vector space $\F_{p^n}$ over $\F_p$;
it is called \textit{primitive}, if one (and hence any) of its roots is a
generator of the multiplicative group $\F_{p^n}^*$.

Here is an example. Consider the extension $\F_{16}:\F_2$ of degree 4
as the quotient ring $\F_2[X]/(X^4+X+1)$. Denoting the class of $X$ in this
quotient by $\a$, we can see that $\a$ is a generator of $\F_{16}$, and the
following table lists the orbits of the Galois group $\Gal(\F_{16}:\F_2)$,
the corresponding irreducible polynomials and, for the polynomials of
maximum degree, indicates their nature:

\begin{center}
\begin{tabular}{|c|c|c|c|}
orbit & polynomial & normal? & primitive? \\
\hline
$\{0\}$ & $x$ && \\
$\{1\}$ & $x+1$ && \\
$\{\a^5,\a^{10}\}$ & $x^2+x+1$ && \\
$\{\a,\a^2,\a^4,\a^8\}$ & $x^4+x+1$ & no & yes \\ 
$\{\a^3,\a^6,\a^{12},\a^9\}$ & $x^4+x^3+x^2+x+1$ & yes & no \\
$\{\a^7,\a^{14},\a^{13},\a^{11}\}$ & $x^4+x^3+1$ & yes & yes
\end{tabular}
\end{center}

\section{Reutenauer's construction}
\label{reut}

By definition, a normal basis is an orbit of the Galois
group which constitutes a basis of $\F_{p^n}$ as a vector space over
$\F_p$. Given a normal basis $A=\{\a,\a^p,a^{p^2},\dots\}$ and a 
necklace $\nu_0,\nu_1,\dots,\nu_{n-1}$, we set up a sum
$\nu_0\a+\nu_1\a^p+\dots+\nu_{n-1}\a^{p^{n-1}}$. Cyclical shifts on the 
sequence
$\nu_0,\nu_1,\dots$ lead to the change of the resulting element of the 
big field
$\F_{p^n}$ within the same orbit of the Galois group, so that the 
minimal polynomial of that element remains the same.

We have therefore two finite sets of equal cardinality, $N_n^p$ and 
$I_n^p$, equipped with a family of 1--1 maps $\f_A:N_n^p\to I_n^p$.
It turns out that the group generated by this set according to the 
construction of Section \ref{gen}, coincides in this case simply with
the set of all maps $\f_A^{-1}\circ\f_B^{}$ of the set of
necklaces into itself. It forms an abelian group, earlier denoted by 
$RG_n^p$, whose order is equal to the
number of normal bases (in the case $p=2$ known as Sloane's sequence
A027362).

Here is the table of automorphism groups of necklaces for $p=2$ and
$n\le15$ ($M_n$ stands for the number of normal bases, that is, the order of 
the group, and the last column lists the lengths of the orbits
of $RG_n^2$ on the corresponding set of necklaces):

$$
\begin{array}{r|r|r|l|l}
 n & |N_n^2| &   M_n  &  \mbox{group}  &     \mbox{orbits}        \\
\hline
 2 &    1  &    1  &   \Z_1                 &    1  \\                     
\hline
 3 &    2  &    1  &   \Z_1                 &    1+1   \\                 
\hline
 4 &    3  &    2  &   \Z_2                 &    2+1      \\              
\hline
 5 &    6  &    3  &   \Z_3                 &    3+3        \\           
\hline
 6 &    9  &    4  &   \Z_2\op\Z_2              &    4+2+2+1    \\           
\hline
 7 &   18  &    7  &   \Z_7                 &    7+7+1+1+1+1  \\         
\hline
 8 &   30  &   16  &   \Z_2\op\Z_2\op\Z_4           &    16+8+4+2     \\         
\hline
 9 &   56  &   21  &   \Z_{21}                &    21+21+7+7    \\         
\hline
10 &   99  &   48  &   \Z_2\op\Z_2\op\Z_2\op\Z_6        &    48+24+24+3   \\         
\hline
11 &   186 &   93  &   \Z_{93}                &    93+93        \\         
\hline
12 &  335  &  128  &   \Z_2\op\Z_2\op\Z_2\op\Z_4\op\Z_4     &    128+64+32+32+16+16\\ 
   &       &       &                      &  +16+8+8+4+4+4+2+1 \\
\hline
13 &  630  &  315  &   \Z_{315}               &    315+315      \\         
\hline
14 &  1161 &  448  &   \Z_2\op\Z_2\op\Z_2\op\Z_2\op\Z_2
                                             & 448+224+224+56+56 \\   
   &       &       &    \op\Z_{14}        &  +28+28+28+28+8+8+7 \\
   &       &       &                      &  +4+4+4+4+1+1 \\
\hline
15 &  2182 &  675  &   \Z_3\op\Z_{15}\op\Z_{15}     &    675+675+225+225 \\       
   &       &       &                      &  +45+45+45+45+45+45 \\
   &       &       &                      &  +15+15+15+15+15+15 \\
   &       &       &                      &  +3+3+3+3+3+3 \\
   &       &       &                      &  +1+1+1+1 \\
\hline
\end{array}
$$

\section{Orbits of the Reutenauer group}
\label{orbits}

It is very interesting to study the \textit{orbits} of the group 
$RG_n^p$ acting on the set of aperiodic necklaces $N_n^p$.
So far, we only have some empirical results in the case $p=2$.

For example,
\begin{itemize}
\item
If $n=4$, the orbits are $O_1=\{0001,0111\}$ and $O_2=\{0011\}$.
The action of the group $RG_4^2=\Z_2$ on each orbit is evident.
\item
If $n=5$, we have $O_1=\{00001,00111,01011\}$ and 
$O_2=\{00011,01111,00101\}$. The action of the group is cyclic on each 
orbit.
\item
If $n=6$, then $O_1=\{000001,011111,001011,001101\}$, 
$O_2=\{000011,010111\}$, $O_3=\{000101,001111\}$ and $O_4=\{000111\}$.
Here, the automorphism group acts as $\Z_2\oplus\Z_2$ on $O_1$, as $\Z_2$ 
on $O_2$ and $O_3$, and trivially on the last orbit.
\end{itemize}

Returning to ideas of Section \ref{intr}, we tried to evaluate various symmetric 
functions of the coefficients in the Drinfeld associator over the orbits 
in these cases, but could not arrive to any sensible conjecture.

It is also worthwhile no notice that there is an interesting operator on 
the set of necklaces, which sometimes turns an aperiodic necklace into a 
periodic one, but in general it takes the whole orbits in the above 
lists into whole orbits. 
We call it the \textit{averaging} operator; by definition, it acts as 
follows: $\{\nu_i\}\mapsto\{\nu_i+\nu_{i+1}\}$. For the above examples, we 
have
\begin{itemize}
\item
$n=4$: $O_1\longrightarrow O_2\longrightarrow\emptyset$.
\item
$n=5$: $O_1\longrightarrow O_2\longrightarrow O_2$.
\item
$n=6$: $O_1\longrightarrow O_2\longrightarrow O_3\longrightarrow O_3$,
$O_4\longrightarrow\emptyset$.
\end{itemize}
(Here, going to $\emptyset$ means that the necklace becomes periodic.)

We noticed that in all examples with $n\le15$, there is a \textit{main} orbit,
that is, an orbit of maximal length equal to the order of the group $RG_n^2$
which acts on this orbit simply transitively.
Iterated averaging operators applied to the main orbit give the majority of orbits; some smallest orbits may go to $\emptyset$.

Remark that here we spoke about the orbits of arithmetical groups
in the sense of Section \ref{grp}. From the point of view of studying the Drinfeld associator, however, it makes little sense to distinguish between symbols 0 and 1: it is more reasonable to extend our groups $RG_n^2$ by the operator that flips 0's and 1's in each necklace. The orbits of such extended groups consist of either one or two orbits of the initial groups;
it should be interesting to have their explicit description and try to study, for example, the sums of coefficients in the logarithm of the associator over these extended orbits.

\section{Sandpile groups}
\label{sand}

Let $\G$ be a finite directed multigraph: it is defined by a finite set of 
vertices $V$ and, for each pair of vertices $v,w\in V$, a non-negative integer
$e(v,w)$, called the number of arrows from $v$ to $w$. 
The total number of
arrows going out of $v$ is referred to as the \textit{outdegree} of $v$ and 
denoted by $\outdeg(v)$; likewise, the \textit{indegree} of $v$ is the
number of arrows going into $v$, denoted by $\indeg(v)$. (Computing these
quantities, we do not take the loops, if any, into consideration.) 
We will assume
that the graph $\G$ is \textit{strongly connected}, that is, there is a directed path
from any vertex $v$ to any other vertex $w$. We will also suppose that our
graph is \textit{Eulerian}, that is, for every vertex $v\in V$ we have
$\indeg(v)=\outdeg(v)$. Under these assumptions, with the given graph $\G$
one can associate a certain finite abelian group $S(\G)$, called the 
\textit{sandpile group} of $\G$ (see \cite{Lev,LP}). The group $S(\G)$ is
defined uniquely up to isomorphism, and the simplest way to define it is
through the \textit{Laplacian matrix} of $\G$.

Let $V=\{v_1,\dots,v_n\}$. The Laplacian matrix $L=(l_{ij})$ of size
$n\times n$ is defined by its entries as
$$
  l_{ii}=-\indeg(v_i),\quad
  l_{ij}=e(v_i,v_j),\mbox{ if } i\ne j.
$$
Let $\L\subset\Z^n$ be the lattice spanned by the rows of $L$.
Evidently, $\L$ is a sublattice of $\Z_0^n=\{(a_1,\dots,a_n)\mid a_1+\dots
+a_n=0\}$. Then we set
$$
  S(\G) = \Z_0^n/\L.
$$

It is known that the group $S(\G)$ can also be
defined as follows. Delete any row and any column from matrix $L$ and call
the resulting $(n-1)\times(n-1)$ matrix $L'$. Let $\L'$ be the sublattice of
$\Z^{n-1}$ spanned by the rows of $L'$. Then $S(\G)\cong\Z^{n-1}/L'$.
On a practical side, to compute the sandpile group, it is enough to reduce
the Laplacian matrix by integral elementary operations on rows and columns 
to its Smith 
normal form, which is a diagonal matrix with integers $(d_1,...,d_{n-1},0)$
on the diagonal; then the group is $\oplus_{i=1}^{n-1}\Z_{d_i}$.

In the next section we will define a series of Eulerian directed
multigraphs $\G_n^p$ labelled by a prime number $p$ and a natural number $n$
and then show that, for $p=2$, their sandpile groups are isomorphic to the 
automorphism groups of necklaces defined above.

\section{Generalized de Bruijn graphs}

Let $\G_n^p$ be the graph with vertex set $V=\Z_n$, the residues modulo $n$,
and $p$ directed edges from every vertex $i$ to each of $pi$, $pi+1$, \dots,
$pi+p-1$. The outdegree of each vertex is thus equal to $p$. It is an easy
exercise to check that the indegree of every vertex is also $p$ and that the
graph is strongly connected. Therefore, the sandpile group $S(\G_n^p)$ is
defined.

We call these graphs \textit{generalized de Bruijn graphs}, because the
well-known de Bruijn graphs appear as a particular case $\G_{2^k}^2$,
see, e.g., \cite{AdB}.

The structure of the group $S(\G_n^2)$ for an arbitrary $n$ is
completely determined by the following two lemmas (Lemma \ref{lemt} 
and Lemma \ref{lemh}), the first of which treats the case of odd $n$ 
and the second shows how to pass from any $n$ to $2n$. Before stating 
the lemmas, let us explain how one could actually arrive at the first,
more difficult, one.
Until the end of this section, we fix $p=2$ and omit the superscripts 2
from various notations.

Suppose that $n$ is odd.
Ideologically, the problem is quite simple: is suffices to find the Smith
normal form of the integer matrix $A_n$ explicitly defined by
$A_n[0,0]=A_n[n-1,n-1]=-1$, $A_n[0,1]=A_n[n-1,n-2]=1$, $A_n[i,i]=-2$ 
for $0<i<n-1$
and $A_n[i,2i]=A_n[i,2i+1]=1$ for $i\in\Z_n$, $i\not=0$, $i\not=n-1$.
The problem is purely
technical, but rather difficult: to understand this, it is enough to look
at the table of first 15 values of the sequence $S(\G_n)$ which coincides
with the table of Section \ref{reut} and is quite non-trivial.
The key difficulty is that, as a rule, there are
three non-zero elements in each column and each row of this matrix, e.g.
$$A_9=\begin{pmatrix}
    -1 &  1 &  0 &  0 &  0 &  0 &  0 &  0 &  0\\
     0 & -2 &  1 &  1 &  0 &  0 &  0 &  0 &  0\\
     0 &  0 & -2 &  0 &  1 &  1 &  0 &  0 &  0\\
     0 &  0 &  0 & -2 &  0 &  0 &  1 &  1 &  0\\
     1 &  0 &  0 &  0 & -2 &  0 &  0 &  0 &  1\\
     0 &  1 &  1 &  0 &  0 & -2 &  0 &  0 &  0\\
     0 &  0 &  0 &  1 &  1 &  0 & -2 &  0 &  0\\
     0 &  0 &  0 &  0 &  0 &  1 &  1 & -2 &  0\\
     0 &  0 &  0 &  0 &  0 &  0 &  0 &  1 & -1
\end{pmatrix}$$
It would be much easier to treat a matrix where there are only two non-zero
elements in each row and column. This goal is almost achieved through a trick
invented by S.H.Chan in his Bachelor's paper \cite{Ch}.

Let us consider the operator given by the Laplace matrix of $\G_n$
in the basis $e_0-e_1$, \dots, $e_{n-2}-e_{n-1}$, $e_{n-1}$, that is, consider the 
matrix $A_n'=C_n^{-1}\cdot A_n\cdot C_n$,
where $C_n$ is a lower triangular matrix of size $n$ with $1$'s on the main 
diagonal and
$-1$'s on the adjacent diagonal. Multiplying a matrix by either $C_n$ or
$C_n^{-1}$ is
equivalent to elementary operations on its rows and columns, hence the
Smith normal forms of matrices $A_n$ and $A_n'$ are the same.
For the previous example we will obtain
$$
A_9'=\begin{pmatrix}
  0 &  0 &  0 &  0 &  0 &  0 &  0 &  0 &  0\\
  1 & -2 &  1 &  0 &  0 &  0 &  0 &  0 &  0\\
  1 &  0 & -2 &  0 &  1 &  0 &  0 &  0 &  0\\
  1 &  0 &  0 & -2 &  0 &  0 &  1 &  0 &  0\\
  1 &  0 &  0 &  0 & -2 &  0 &  0 &  0 &  1\\
  0 &  1 &  0 &  0 &  0 & -2 &  0 &  0 &  0\\
  0 &  0 &  0 &  1 &  0 &  0 & -2 &  0 &  0\\
  0 &  0 &  0 &  0 &  0 &  1 &  0 & -2 &  0\\
  0 &  0 &  0 &  0 &  0 &  0 &  0 &  1 & -2
\end{pmatrix}$$
We see that the lower-right minor of codimension 1 has the required
property, and its Smith normal form can be found by drawing horizontal and vertical lines
between the non-zero entries in each row and column of the matrix and
considering the cycles obtained:
\newcommand{\cdt}{\cdots}
\newcommand{\vdt}{\vdots}
$$\begin{array}{ccccccccccccccc}
 -2 &\cdt&  1 &    &  0 &    &  0 &    &  0 &    &  0 &    &  0 &    &  0 \\
\vdt&    &\vdt&    &    &    &    &    &    &    &    &    &    &    &    \\
  0 &    & -2 &\cdt&  0 &\cdt&  1 &    &  0 &    &  0 &    &  0 &    &  0 \\
\vdt&    &    &    &    &    &\vdt&    &    &    &    &    &    &    &    \\
  0 &    &  0 &    & -2 &\cdt&  0 &\cdt&  0 &\cdt&  1 &    &  0 &    &  0 \\
\vdt&    &    &    &\vdt&    &\vdt&    &    &    &\vdt&    &    &    &    \\
  0 &    &  0 &    &  0 &    & -2 &\cdt&  0 &\cdt&  0 &\cdt&  0 &\cdt&  1 \\
\vdt&    &    &    &\vdt&    &    &    &    &    &\vdt&    &    &    &\vdt\\
  1 &\cdt&  0 &\cdt&  0 &\cdt&  0 &\cdt& -2 &    &  0 &    &  0 &    &  0 \\
    &    &    &    &\vdt&    &    &    &\vdt&    &\vdt&    &    &    &\vdt\\
  0 &    &  0 &    &  1 &\cdt&  0 &\cdt&  0 &\cdt& -2 &    &  0 &    &  0 \\
    &    &    &    &    &    &    &    &\vdt&    &    &    &    &    &\vdt\\
  0 &    &  0 &    &  0 &    &  0 &    &  1 &\cdt&  0 &\cdt& -2 &    &  0 \\
    &    &    &    &    &    &    &    &    &    &    &    &\vdt&    &\vdt\\
  0 &    &  0 &    &  0 &    &  0 &    &  0 &    &  0 &    &  1 &\cdt& -2
\end{array}$$
In this example we see two cycles of lengths 4 and 12 which are simply a
visualization of the orbits of lengths 2 and 6 in the set
$\Z_9\setminus\{0\}$ under the doubling operator $x\mapsto 2x$ (in our case
the orbits are $\{1,2,4,8,7,5\}$ and $\{3,6\}$). It is readily verified that
each orbit of length $d$ adds a summand $\Z_{2^d-1}$ to the Smith group
of such a matrix (where each row and each column contain one entry 1 and 
one entry $-2$), so for our example we obtain $\Z_{63}\oplus\Z_3$.
Unfortunately, the presence of a nonzero first column spoils this clear
picture, namely, it decreases the size of the group by a factor of $n$
(more exactly, it leads to a subgroup of index $n$). For the example under
study, any subgroup of index 9 is isomorphic to $\Z_{21}$. However, there
are situations where such a group may have different subgroups of index $n$.
S.H.Chan in a series of rather involved technical lemmas showed how exactly
looks the resulting group $S(\G_n)$. In most cases, one must simply divide
by $n$ the order of the first cyclic group, corresponding to the orbit of
the number 1. This is so for all odd integers up to 19. For $n=21$, however,
the set $\Z_{21}\setminus\{0\}$ decomposes into five orbits
$\{1,2,4,8,16,11\}$, $\{3,6,12\}$, $\{5,10,20,19,17,13\}$, $\{7,14\}$ and
$\{9,18,15\}$ of lengths 6, 3, 6, 2 and 3, respectively, but the sandpile
group is actually equal to $\Z_9^2\oplus\Z_7^3$, and not to
$\Z_3^2\oplus\Z_9\oplus\Z_7^3$ as one might infer from the previous rule.
To state the exact formula proven by S.H.Chan, we need some notations.

Let $n$ be an odd number. For any element $v\in\Z_n\setminus0$ let 
$l(v)$ be the length of its orbit under the doubling operator 
$x\mapsto2x$. Let $H_n$ be the set of minimal representatives of all 
orbits. Now, denote by $\P_n$ the set of all prime divisors of $n$. 
For each $q\in\P_n$ let $q'$ stand for the maximal power of $q$ that 
divides $n$ and let $q''=n/q'$. Denote the set of all such 
residues $q''$ by $V_n$. Finally, for an abelian group $G$ and an 
integer $k$ let $kG$ be the subgroup of all elements of $G$ of the form 
$kx$, $x\in G$. 
Then

\begin{lemma}
\label{lemt}
If $n$ is odd, then the group $S(\G_n)$ has the following
decomposition
$$
  S(\G_n) \cong 
\bigoplus_{q \in \P_n } q' Z_{2^{l(q'')}-1}
\oplus
\bigoplus_{v \in H_n \setminus V_n} Z_{2^{l(v)}-1}
$$
\end{lemma}

Example. Take $n=21$.
Then $\P_n=\{3,7\}$,
$q'_1=3$, $q''_1=7$, $l(7)=2$,
$q'_2=7$, $q''_2=3$, $l(3)=3$,
$H_{21}=\{1,3,5,7,9\}$, $V_{21}=\{3,7\}$, so the first direct summand in the
above formula is trivial, and the second gives:
$$
  S(\G_{21})= \Z_{63}\oplus\Z_{63}\oplus\Z_7.
$$

The proof of the next lemma is much simpler: it follows from the fact 
that $\G_{2n}$ is the directed line graph of $\G_n$ (see \cite{Lev}).

\begin{lemma}
Suppose that $n=2^km$ where $m$ is odd. Then
$$
  S(\G_n) \cong S(\G_m) 
  \oplus \Z_{2^k}^{m-1} \oplus 
  \left[ \bigoplus_{i=2}^{k} \mathbb{Z}_{2^{k+1-i}}^{2^{i-2}m}
  \right]
$$\label{lemh}
\end{lemma}

\section{Circulant matrices}
\label{circ}
An $n\times n$ matrix $A=a_{i,j}$, $i,j\in\Z_n$, over a field  $\K$ 
is called {\em circulant} if its
rows are the cyclic shifts of the first row, i.e. 
$$ 
A_{i+1, j+1}=A_{ij},\qquad \text{for all}\ i,j\in\Z_n,
$$
where $i+1$ and $j+1$ are taken modulo $n$.

In particular, the permutation matrices associated with the powers
of the cyclic permutation $(0,1,\dots,n-1)$ are
circulant; they form a basis of the algebra $C_n(\K)\cong\K[\Z_n]$
of all $n\times n$ circulant matrices over $\K$. We see that
the algebra $C_n(\K)$ has dimension $n$ and is commutative.

For a circulant matrix to be non-degenerate it is necessary (but not
sufficient) that its rows (and columns) are aperiodic. Denote the group
of non-degenerate circulants by $C_n(\K)^*$. By the observation made above, it is commutative.
In the case $\K=\F_p$, by studying the natural action of this group on the 
field $\F_{p^n}$ considered as a vector space over $\F_p$, it is easy to
deduce that the Reutenauer group $RG_n^p$ is isomorphic to the quotient 
$C_n(\F_p)^*/\Z_n$,
where $\Z_n$ is the group of permutation circulant matrices, those
associated with the powers of $(0,1,\dots,n-1)$.

In the case $p=2$, it was proved in \cite{Ch} that the series of these 
quotient groups satisfies the same relations as those given in lemmas 
\ref{lemt} and \ref{lemh}. The proof of these facts is not so involved 
as the proof of lemma \ref{lemt} and relies basically on the primary 
decomposition theorem from linear algebra. The main theorem follows.

\begin{thm} (S.W.Chan) For any natural $n$ we have
$$
   RG_n^2 \cong S(\G_n^2)\,.
$$
\end{thm}

This theorem is quite remarkable, because it relates the objects coming 
from entirely different areas of mathematics. It is noteworthy that nobody 
knows any explicit isomorphism between the two groups in question, although
the elements of both can be encoded by some sequences of 0's and 1's.

In the paper \cite{CHP}, this result is generalized to any prime 
number $p$ as follows: $RG_n^p \cong S(\G_n^p)\oplus\Z_{p-1}$.
Moreover, that paper describes the 
structure of sandpile groups for the generalized de Bruijn graphs $\G_n^p$
for arbitrary values of $p$, not only prime. Of course, in the general case
there is no analogue of the isomorphism theorem; 
however, if $p$ is a power of a prime both groups are defined
and the relation between them should be studied.

\section{Golomb's construction}
\label{gol}

We conclude the paper with some experimental data related to another
set of 1--1 correspondences between the necklaces and irreducible polynomials
mentioned above. 

Let $\a$ be a generator of the multiplicative group of
the field $\F_q$, $q=p^n$. 
S.\,Golomb \cite{Gol} defined a bijection
$\psi_\a:N_n^p\to I_n^p$ which depends only on the orbit of the element $\a$
under the Galois group action, that is, $\psi_\a$ is completely determined by
the primitive polynomial with one of the roots $\a$.
To a necklace $\nu=(\nu_0,...,\nu_{n-1})$ we assign the following element of
$\F_q$:
$$
\a^{\nu_0+p\nu_1+...+p^{n-1}\nu_{n-1}}
$$ 
and then take its minimal polynomial, which we denote by $\psi_\a(\nu)$.
It is not hard to prove (see \cite{Gol}) that the map $\psi_\a$ is 1-to-1
for any $\a$ which is a root of a primitive polynomial, and that these maps
are the same for all roots of one primitive polynomial, and are distinct for different primitive polynomials, so that they generate a group
of automorphisms of necklaces whose order equals the number of primitive
polynomials $\phi(2^n-1)/n$, where $\phi$ is Euler's totient function.
This group turns out to be abelian, too.

Here is a table of these groups for $p=2$ and $2\le n\le12$ given 
together with the sizes of the orbits into which they split the set of
necklaces ($M_n$ stands for the order of the group):

$$\begin{array}{r|r|r|l|l}
 n & |N_n^2| &   M_n  &  \mbox{group}  &     \mbox{orbits}        \\
\hline
 2  &  1 &   1 & \Z_1             &  1\\
\hline
 3  &  2 &   2 & \Z_2             &  2\\
\hline
 4  &  3 &   2 & \Z_2             &  2+1\\
\hline
 5  &  6 &   6 & \Z_6             &  6\\
\hline
 6  &  9 &   6 & \Z_6             &  6+2+1\\
\hline
 7  & 18 &  18 & \Z_{18}          &  18\\
\hline
 8  & 30 &  16 & \Z_2\oplus\Z_8        &  16+8+4+2\\
\hline
 9  & 56 &  48 & \Z_2\oplus\Z_{24}     &  48+8\\
\hline
10  & 99 &  60 & \Z_2\oplus\Z_{30}     &  60+30+6+2+1\\
\hline
11  &186 & 176 & \Z_2\oplus\Z_{88}     &  176+8+2\\
\hline
12  &335 & 144 & \Z_2\oplus\Z_6\oplus\Z_{12}&  144+48+36+2\cdot24+2\cdot12\\
    &    &     &                  & +8+2\cdot6+2\cdot4+3\cdot2+1\\
\end{array}$$

Notice that the sequence of orders of these groups is not monotonous.
As yet, nobody knows any relations between these groups and other 
mathematical objects, which was the case for the Reutenauer groups.

\section{Acknowledgments}

We are grateful to M.\,Vsemirnov who draw our attention to the close relation between the Reutenauer groups and circulant matrices.
The first author acknowledges the hospitality of the Nanyang Technological University, Singapore, and Max-Planck Institute for Mathematics, Germany, 
where this paper was accomplished.

\vspace{1cm}

\begin{flushleft}
St.Petersburg Division of the\\
Steklov Mathematical Institute\\
e-mail \texttt{duzhin@pdmi.ras.ru}
\vspace{5mm}

Nanyang Technological University, Singapore\\
e-mail \texttt{dima@ntu.edu.sg}
\end{flushleft}

\end{document}